\documentclass[9pt,technote]{IEEEtran}
\usepackage{graphicx}          
\usepackage{amsmath}           
\usepackage{amssymb}
\usepackage{algorithm}
\usepackage{algpseudocode}
\usepackage[latin1]{inputenc}  

\usepackage{amsthm}

\theoremstyle{definition}

\DeclareMathOperator{\E}{E}
\DeclareMathOperator*{\argmin}{arg\,min}
\DeclareMathOperator*{\argmax}{arg\,max}

\begin{document}

\title{Dynamic Iterative Pursuit}
\author{Dave Zachariah, Saikat Chatterjee and Magnus Jansson\thanks{
The authors are with the ACCESS Linnaeus Centre, KTH-Royal Institute of
Technology, Stockholm. E-mail:
$\{$dave.zachariah,magnus.jansson$\}$@ee.kth.se and
saikatchatt@gmail.com. This work was partially supported by the Swedish Research Council under contract 621-2011-5847.}}

\maketitle
\begin{abstract}
For compressive sensing of dynamic sparse signals, we develop an iterative pursuit algorithm. A dynamic sparse signal process is characterized by varying sparsity patterns over time/space. For such signals, the developed algorithm is able to incorporate sequential predictions, thereby providing better compressive sensing recovery performance, but not at the cost of high complexity. Through experimental evaluations, we observe that the new algorithm exhibits a graceful degradation at deteriorating signal conditions while capable of yielding substantial performance gains as conditions improve.
\end{abstract}

\section{Introduction}

Compressive Sensing (CS) \cite{Candes&Wakin2008} problems assume a
sparse-signal model, undersampled by a linear measurement process.
The algorithms for CS can be separated into three broad classes:
convex relaxation, Bayesian inference, and iterative pursuit
(IP). For large-dimensional CS signal-reconstruction, IP
algorithms offer computationally efficient solutions. Examples of
such IP algorithms are orthogonal matching pursuit (OMP)
\cite{Tropp&Gilbert2007}, subspace pursuit (SP)
\cite{Dai&Milenkovic2009} and several variants of them
\cite{DonohoEtAl2006}\cite{Needell&Tropp2009}\cite{Needell&Vershynin2009}\cite{ChatterjeeEtAl2011b}.
The methodology of such IP algorithms is to detect and
reconstruct the non-zero, or `active', signal coefficients in a least-squares
framework. These algorithms may use some prior information, such
as the maximum allowable cardinality of the `support set'. The
support set is defined as the set of active signal coordinates
of the underlying sparse signal. In general, the IP algorithms
work with a single snapshot of the measurements. In this paper, we
are interested in generalizing the iterative pursuit approach so as to use
more prior information. Such prior information may for instance be
available in dynamically evolving sparse processes with
temporal/spectral/spatial correlations, as in the sparse signal
scenarios of magnetic resonance imaging (MRI) \cite{LustigEtAl2007,
  Vaswani2010}, spectrum sensing \cite{SundmanEtAl2010} and direction of arrival estimation \cite{MalioutovEtAl2005}.

Incorporation of prior information is a recent trend in CS. In
\cite{Vaswani2008}, the overall methodology is sequential and can
be seen as a two-step approach: (1) support-set detection of the
sparse signal, and (2) reduced-order recovery using prior
information on the detected support set. For a reasonable
detection of support set, \cite{Vaswani2008} uses convex
relaxation algorithms. Then, a standard Kalman filter (KF) is
employed to use prior information for sequential signal recovery.
Without explicit support set detection, \cite{CarmiEtAl2010} uses
KF to estimate the entire signal and enforces sparsity by imposing
an approximate norm constraint. However, the work of
\cite{CarmiEtAl2010} validates their algorithm for a signal with a
static sparsity pattern (i.e. an unknown pattern that does not
evolve over time). Similarly, \cite{Zhang&Rao2010} considers
scenarios with static sparsity patterns and solves the
reconstruction of a temporally evolving sparse signal with
multiple measurement vectors in a batch Bayesian learning
framework with unknown model parameters. Iterative pursuit algorithms that
can use prior information to recover dynamic sparse signals are,
however, largely unexplored. One exception is \cite{DaiEtAl2011} which
uses a maximum aposteriori criterion to modify SP for Gaussian
processes. Their signal model, however, does not allow explicit
modeling of the temporal correlation of the sparsity pattern.

In this paper, we consider a signal model with dynamically evolving
sparsity pattern. In other words, we consider that the signal sparsity pattern varies over time/space at any rate (i.e. from a slowly varying case to a rapidly varying case). We then develop a predictive orthogonal matching
pursuit algorithm that can incorporate prior information in a
stochastic framework, using the signal to prediction error ratio as a
statistic. Thereby recovery performance can be improved while
maintaining the complexity advantage of IP algorithms. This generalizes the
linear minimum mean square error (MMSE) approach taken in
Gaussian-based matching pursuit \cite{ChatterjeeEtAl2011a}.

We also develop a robust detection strategy for finding the support set
elements of a dynamic sparse signal. Compared to standard
correlation-based successive detection in existing iterative pursuit
algorithms, such as OMP, this detection strategy is found to be more
robust to erratic changes in the sparsity pattern.

Finally, the algorithm is integrated into a recursive Kalman-filter
framework in which the sparse process is predicted as a superposition
of state transitions. The new IP algorithm using sequential predictions
is referred to as dynamic iterative pursuit (DIP). Through
experimental simulations we show that the new algorithm provides a
graceful degradation at higher measurement noise levels and/or lower
measurement signal dimensions, while capable of yielding substantial
gains at more favorable signal conditions.

\emph{Notation:} $\| \mathbf{x} \|_0$ denotes $l_0$ `norm', i.e. the
number of non-zero coefficients of the vector $\mathbf{x}$.
$\mathbf{A} \oplus \mathbf{B}$ is the direct sum of matrices.
$|S|$ and $S^c$ are the cardinality and complement of set $S$,
respectively. $\varnothing$ denotes the empty set. $(\cdot)^*$ is the Hermitian transpose operator. $\mathbf{A}^\dagger$ the Moore-Penrose pseudoinverse of matrix
$\mathbf{A}$. $\mathbf{C}^{1/2}$ denotes a matrix square root of a
positive definite matrix $\mathbf{C}$, and $\mathbf{C}^{*/2}$ is
its Hermitian transpose. $\mathbf{A}_{[\mathcal{I},\mathcal{J}]}$
denotes a submatrix of $\mathbf{A}$ with elements from row and
column indices listed in ordered sets $\mathcal{I}$ and
$\mathcal{J}$. Similarly, the column vector
$\mathbf{x}_{[\mathcal{I}]}$ contains the elements of $\mathbf{x}$
with indices from set $\mathcal{I}$.

\section{Signal model}

We consider a standard CS measurement setup,
\begin{equation}
\mathbf{y}_t = \mathbf{H} \mathbf{x}_t + \mathbf{n}_t \in
\mathbb{C}^M,
\end{equation}
where $\mathbf{x}_t \in \mathbb{C}^N$ is the sparse state
vector to be estimated and the $\mathbf{n}_t$ is zero-mean Gaussian,
$\E[\mathbf{n}_t \mathbf{n}^*_{t-l}] = \mathbf{R}_t \delta(l)$. The
sensing matrix $\mathbf{H} = [\mathbf{h}_1 \quad \cdots \quad
\mathbf{h}_N] \in \mathbb{C}^{M \times N}$, where $M < N$. Both $\mathbf{H}$ and $\mathbf{R}_t$ are given. Without loss of generality, we assume $\| \mathbf{h}_i \|_{2}  = 1$.
\subsection{Process model}

Let the `support set' $I_{x,t} \subset \{1, \dots, N\}$ represent the
sparsity pattern of $\mathbf{x}_t \in \mathbb{C}^N$. It will be assumed
that $ |I_{x,t}| \equiv \| \mathbf{x}_t \|_0 \leq K_{\text{max}}$,
where $K_{\text{max}} < M$. Let $\lambda_{ji}$ denote the state
transition probability $j \rightarrow i$ of the `active' signal
coordinate $j$. Then the probabilities determine the transition
$I_{x,t} \rightarrow I_{x,t+1}$, as will be illustrated below.

The transition of an active signal coordinate $j
\rightarrow i$ is modeled as an autoregressive (AR) process,
\begin{equation}
x_{i,t+1} = \alpha_{ij} x_{j,t} + w_{i,t},
\label{eq:linearjump}
\end{equation}
where $x_{j,t}$ denotes the $j$th component of $\mathbf{x}_t$,
$w_{i,t}$ is the associated innovation and the AR coefficient $|
\alpha_{ij}| < 1$. This model extends the scenario considered in
\cite{Zhang&Rao2010} where the transition probabilities are degenerate
$\lambda_{ji} = \delta(i-j)$, resulting in a static sparsity pattern
$I_{x,t} \equiv I_x, \forall t$.

The sparse-signal process can be written
compactly as a linear state-space model with random transition
matrices $\mathbf{A}_t$ and $\mathbf{B}_t$,
\begin{equation}
\mathbf{x}_{t+1} = \mathbf{A}_t \mathbf{x}_t + \mathbf{B}_t
\mathbf{w}_t, \label{eq:processmodel}
\end{equation}
where $\mathbf{w}_t$ is zero-mean Gaussian, $\E[\mathbf{w}_t
\mathbf{w}^*_{t-l}] = \mathbf{Q}\delta(l) \in \mathbb{C}^{N
  \times N}$ and $\mathbf{Q} = \text{diag}(\sigma^2_1, \dots,
\sigma^2_N)$. The non-zero elements of $\mathbf{A}_t \in \mathbb{C}^{N
  \times N}$ are $a_{ij,t} = \alpha_{ij}$ for all $j \in I_{x,t}$ and
$i \in I_{x,t+1}$. Similarly, the non-zero elements of the diagonal
matrix $\mathbf{B}_t \in \mathbb{C}^{N \times N}$  are $b_{ii,t} = 1$
for all $i \in I_{x,t+1}$. The model parameters $\alpha_{ij},
\lambda_{ji}$ and $\mathbf{Q}$ are assumed to be known.

\subsection{Examples}

Use of the transition probabilities $\lambda_{ji}$ along with signal model \eqref{eq:linearjump} enables compact modeling of dynamically evolving
sparsity patterns. The potential applications include MRI, spectrum sensing, direction of arrival estimation, frequency tracking etc.

\begin{figure}
  \begin{center}
    \includegraphics[width=0.90\columnwidth]{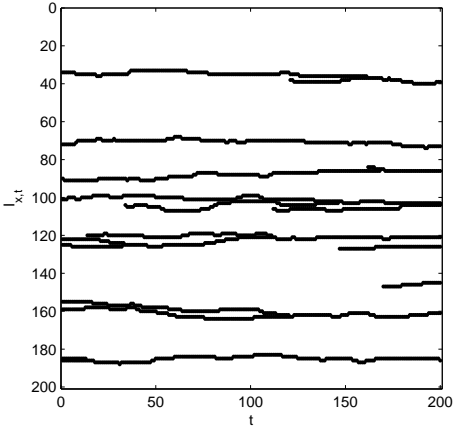}
  \end{center}
  \vspace*{-0.6cm}
  \caption{Example of evolving sparsity pattern with $(N,K,T) = (200,10,200)$ and transition probabilities \eqref{eq:lambda_model_1}.}
  \label{fig:examplepattern_1}
\end{figure}

As an initial example, consider a slowly varying sparsity pattern
$I_{x,t}$ over $T$ snapshots, following\footnote{The two cases $j = 1$
  or $j=N$ are necessary for the edge states.}
\begin{equation}
\lambda_{ji} =
\begin{cases}
0.90 & i=j \\
0.05 & i=j\pm1,\text{ if } j \not \in \{1,N\} \\
0.10 & i=j+1,\text{ if } j = 1 \\
0.10 & i=j-1,\text{ if } j = N \\
0      & |i-j| > 1
\end{cases}.
\label{eq:lambda_model_1}
\end{equation}
A realization of this process is illustrated in
Figure~\ref{fig:examplepattern_1}. This choice is intended to model the
strong temporal correlation of sparse signals exhibited in e.g.
MRI.

\begin{figure}
  \begin{center}
    \includegraphics[width=0.95\columnwidth]{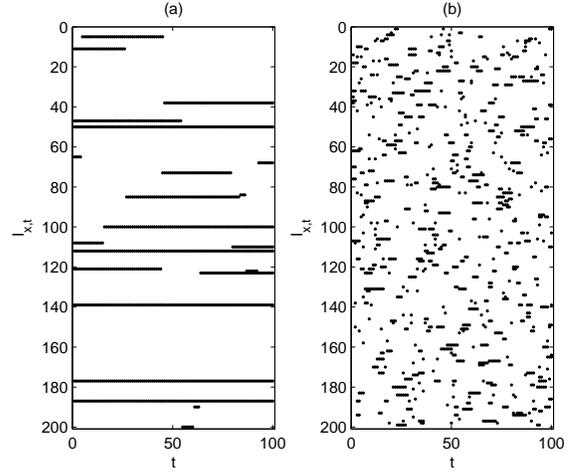}
  \end{center}
  \vspace*{-0.6cm}
  \caption{Examples of evolving sparsity patterns with $(N,K,T) = (200,10,100)$ and transition probabilities \eqref{eq:lambda_model_2} with (a) $\nu = 0.01$ and (b) $\nu = 0.5$.}
  \label{fig:examplepattern_2}
\end{figure}

Next, consider a simpler parameterization,
\begin{equation}
\lambda_{ji} =
\begin{cases}
1 - \frac{N-1}{N}\nu & i=j \\
\frac{1}{N} \nu & i \neq j
\end{cases},
\label{eq:lambda_model_2}
\end{equation}
where $\nu \in [0,1]$ is a mixture factor. This is intended to model
more erratically evolving patterns in e.g. frequency-hopping radio
frequency (RF) signals. Examples
of resulting sparsity patterns are shown in Figure
\ref{fig:examplepattern_2}, where we consider $\nu = 0.01$ and 0.5. It
can be seen that the evolution of the sparsity pattern becomes more
erratic as $\nu$ increases. In the above examples we have ensured that the sparsity level is constant, $K=10$.

\section{Dynamic Iterative Pursuit}

We approach the dynamic estimation problem by first developing an
iterative pursuit algorithm that can incorporate prior information in
the form of a prediction of $\mathbf{x}_t$. A support set detection
strategy is proposed using the signal to prediction error ratio. Next
we develop a recursive algorithm based on the Kalman-filter
framework. We propose predicting the sparse process as the
superposition of all state transitions of the signal coefficients.

\subsection{Incorporation of prior information}


Given the constraint on the support set, $| I_x | \leq K_{\text{max}}$, the
brute force least-squares solution would be to enumerate all
combinations of possible support sets. For each set, $I \subset \{ 1,
\dots, N \}$, the signal coefficients are reconstructed by a least-squares
criterion and a measurement residual is computed, $\mathbf{r} =
\mathbf{y} -\mathbf{H}_{[\cdot,I]} \hat{\mathbf{x}}_{[I]}$. The
reconstruction with minimum residual norm is then chosen as the
solution. However, with at least one active coefficient there are
$\sum^{K_{\text{max}} }_{K=1} \binom{N}{K}$ possible support sets $I$
to enumerate, which is clearly intractable.

Several iterative pursuit algorithms solve the estimation problem by a
sequential detection of the support set and reconstruction of the
corresponding signal coefficients. We will use OMP to illustrate the
essential components of this sequential strategy.

OMP takes a support set $I$ as its starting
point. Reconstructed signal coefficients, $\hat{x}_j$, $j \in I$, are
cancelled from the observation $\mathbf{y}$ to form the residual
$\mathbf{r} = \mathbf{y} -\mathbf{H}_{[\cdot,I]}
\hat{\mathbf{x}}_{[I]}$. Initially $I =\varnothing$. Under the
hypothesis of a remaining active coefficient $x_i$, $i \not \in I$,
the residual signal model is
\begin{equation}
\begin{split}
\mathbf{r} &=  \mathbf{h}_i x_i + \sum_{j \in I}   \mathbf{h}_j \xi_j +  \mathbf{n},
\end{split}
\label{eq:residualhyp}
\end{equation}
where $\xi_j = x_{j} - \hat{x}_{j}$ are estimation errors. OMP detects the active coefficient by using a matched filter. The matched filter employs the strategy of estimating $x_i$ using a least-squares criterion, $\check{x}_i = \mathbf{h}^\dagger_i \mathbf{r} = \mathbf{h}^*_i \mathbf{r}$.
The index $i \not \in I$ corresponding to maximum energy $| \check{x}_i |^2$
is added to $I$. Finally the coefficients corresponding to $I$ are
estimated jointly based on a least-squares criterion, solving
\begin{equation*}
\hat{\mathbf{x}}_{[I]} = \argmin_{\mathbf{x}_{[I]}\in
\mathbb{C}^{|I|}} \left\| \mathbf{y} - \mathbf{H}_{[\cdot,I]}
\mathbf{x}_{[I]} \right\|^2_2 = \mathbf{H}^\dagger_{[\cdot, I]}
    \mathbf{y}.
\end{equation*}
The residual $\mathbf{r}$ is updated and the process is repeated until
the residual norm no longer decreases or when $|I|$ reaches the limit
$K_{\text{max}}$. For sake of clarity OMP is summarized in Algorithm~\ref{alg:OMP} where $k$ denotes the iteration index.
\begin{algorithm}
\caption{: Orthogonal Matching Pursuit (OMP)} \label{alg:OMP}
\begin{algorithmic}[1]
    \State Given: $\mathbf{y}$ and $\mathbf{H}$
    \State Set $k=0$, $\mathbf{r}_0 = \mathbf{y}$ and $I = \varnothing$
    \Repeat
    \State $k := k+1$
    \State $i_k = \argmax_{i \in I^c}  | \mathbf{h}^*_i
    \mathbf{r}_{k-1} |$
    \State $I := I \cup i_k$
    \State $\hat{\mathbf{x}}_{[I]} = \mathbf{H}^\dagger_{[\cdot, I]}
    \mathbf{y}$; $\hat{\mathbf{x}}_{[I^c]} = \mathbf{0}$
    \State $\mathbf{r}_{k} = \mathbf{y} - \mathbf{H}_{[\cdot,I]} \hat{\mathbf{x}}_{[I]}$
    \Until{$(\| \mathbf{r}_{k} \|_2 \geq \| \mathbf{r}_{k-1} \|_2)$ or $(k
      > K_{\text{max}})$}
    \State Output: $\hat{\mathbf{x}}$ and $I$
\end{algorithmic}
\end{algorithm}

Using a stochastic framework, we now extend the estimation strategy to a
scenario in which a prediction $\hat{\mathbf{x}}^- = \mathbf{x} +
\mathbf{e}$ is given, where $\mathbf{e} \sim \mathcal{N}(\mathbf{0},
\mathbf{P}^-)$ and error covariance matrix $\mathbf{P}^-$ is
known. Then the signal to prediction error ratio,
\begin{equation}
\rho_i \triangleq \frac{ \E[ | x_i |^2 ] }{  \E[ | e_i |^2 ] } \in [0,
\infty),
\end{equation}
quantifies the certainty that $i$ belongs to the support set. We propose to use
$\rho_i$ for selecting indices to be added to $I$. The ratio $\rho_i$ is
successively updated by conditioning the expectations on the residual,
under the hypothesis with signal model \eqref{eq:residualhyp}. Then $\E[ |
x_i |^2 ] = |\mu_{i|r}|^2 + \sigma^2_{i|r}$ where the conditional mean
$\mu_{i|r}$ is given by the MMSE-estimator and $\sigma^2_{i|r}$ by its error
variance. The prior of $x_i$ is the prediction $\hat{x}^-_i$.

For tractability the estimation errors $\xi_j$ are assumed to be
Gaussian and their correlations negligible so that the MMSE-estimator gives
\begin{equation}
\begin{split}
\mu_{i|r} &\simeq \hat{x}^-_i + \mathbf{g}^*_i \left(  \mathbf{r} - \mathbf{h}_i \hat{x}^-_i \right) \\
\sigma^2_{i|r}     &\simeq \left( 1 - \mathbf{g}^*_i \mathbf{h}_i \right) p^-_i,
\end{split}
\label{eq:SPERapprox}
\end{equation}
where $p^-_i$ is the $i$th diagonal element of $\mathbf{P}^-$
\cite{KailathEtAl2000}. The gain (row) vector $\mathbf{g}^*_i$ and covariance matrix $\mathbf{D}$ are,
\begin{equation}
\begin{split}
\mathbf{g}^*_i &= \left (  \frac{1}{p^-_i} + \mathbf{h}^*_i
  \mathbf{D}^{-1} \mathbf{h}_i \right)^{-1} \mathbf{h}^*_i
\mathbf{D}^{-1}, \\
\mathbf{D} &= \sum_{j \in I} \sigma^2_j \mathbf{h}_j \mathbf{h}^*_j +
\mathbf{R},
\end{split}
\end{equation}
where $\sigma^2_j$ is the variance of $\xi_j$ and $\mathbf{R}$
is the covariance matrix of $\mathbf{n}$. As the support set $I$
successively grows, the inverse $\mathbf{D}^{-1}$ can be updated
efficiently using the Sherman-Morrison formula, as shown below.

To sum up, the signal to prediction error ratio is given by
\begin{equation}
\rho_i = \frac{ |\mu_{i|r}|^2 + \sigma^2_{i|r} }{  p^-_i }
\label{eq:SPER}
\end{equation}
and approximated using \eqref{eq:SPERapprox}. For the maximum $\rho_i$, $i
\not \in I$ is added to $I$. Finally, signal coefficients are jointly
re-estimated, solving a weighted least-squares problem
\begin{equation} \hat{\mathbf{x}}_{[I]} = \argmin_{\mathbf{x}_{[I]}\in
\mathbb{C}^{|I|}} \left\|
\begin{bmatrix} \mathbf{y} \\ \hat{\mathbf{x}}^-_{[I]}
\end{bmatrix} -
\begin{bmatrix}
\mathbf{H}_{[\cdot, I]} \\ \mathbf{I}_{|I|} \end{bmatrix}
\mathbf{x}_{[I]} \right\|^2_{\mathbf{R}^{-1} \oplus
\mathbf{S}^{-1} }, \label{eq:MMSE}
\end{equation}
where $\mathbf{S} = \mathbf{P}^-_{[I,I]}$. This is the linear MMSE
estimator provided $I$ is the correct support set. The residual
$\mathbf{r}$ is updated and the process is repeated as above. The
resulting algorithm is referred to as `Predictive OMP' (PrOMP) and is
summarized in Algorithm \ref{alg:PrOMP}.

\begin{algorithm}
\caption{: Predictive Orthogonal Matching Pursuit (PrOMP)}
\label{alg:PrOMP}
\begin{algorithmic}[1]
 \State Given: $\mathbf{y}, \mathbf{H}, \mathbf{R}^{-1},
 \hat{\mathbf{x}}^-$ and $\mathbf{P}^-$
    \State Set $k=0$, $\mathbf{r}_0 = \mathbf{y}$, $I = \varnothing$
    and $\mathbf{D}^{-1} = \mathbf{R}^{-1}$
    \Repeat
    \State $k := k+1$
    \State Compute $\rho_i$ using \eqref{eq:SPER} and \eqref{eq:SPERapprox}
    \State $i_k = \argmax_{i \in I^c}  \rho_i$
    \State $I := I \cup i_k$
   \State $\hat{\mathbf{x}}_{[I]}
   =$\texttt{mmse-rec}$(\mathbf{y},\mathbf{H},\mathbf{R}^{-1},\hat{\mathbf{x}}^-,\mathbf{P}^-,I)$;
   $\hat{\mathbf{x}}_{[I^c]} = \mathbf{0}$
    \State $\mathbf{r}_{k} = \mathbf{y} - \mathbf{H}_{[\cdot,I]}
    \hat{\mathbf{x}}_{[I]}$
   \State $\mathbf{D}^{-1} =$\texttt{update-cov}($\mathbf{D}^{-1},
   \mathbf{H},\mathbf{R}^{-1}, \mathbf{P}^-,I, i_k$)
    \Until{$(\| \mathbf{r}_{k} \|_2 \geq \| \mathbf{r}_{k-1} \|_2)$ or $(k
      > K_{\text{max}})$}
    \State Output: $\hat{\mathbf{x}}$ and $I$
\end{algorithmic}
\end{algorithm}

The function \texttt{mmse-rec} solves \eqref{eq:MMSE} and can be
computed by a measurement update of form:
\begin{equation*}
\hat{\mathbf{x}}_{[I]} = \hat{\mathbf{x}}^-_{[I]} +
\mathbf{K} \left( \mathbf{y} - \mathbf{H}_{[\cdot,I]}
\hat{\mathbf{x}}^-_{[I]} \right),
\end{equation*}
where $\mathbf{K} = ( \mathbf{S}^{-1} + \mathbf{H}^*_{[\cdot,I]}
\mathbf{R}^{-1} \mathbf{H}_{[\cdot,I]}  )^{-1}
\mathbf{H}^*_{[\cdot,I]}  \mathbf{R}^{-1}$. The function
\texttt{update-cov} updates the inverse covariance matrix for the
added reconstructed coefficient $\hat{x}_i$,
\begin{equation*}
\mathbf{D}^{-1} := \mathbf{D}^{-1} \left( \mathbf{I}_M - \frac{\mathbf{h}_i
  \mathbf{h}^*_i \mathbf{D}^{-1}}{ \sigma^{-2}_{i} + \mathbf{h}^*_i
  \mathbf{D}^{-1} \mathbf{h}_i } \right),
\end{equation*}
where $\sigma^{2}_{i}$ is the corresponding diagonal element of the
posterior error covariance matrix $(
\mathbf{S}^{-1} + \mathbf{H}^*_{[\cdot,I]} \mathbf{R}^{-1}
\mathbf{H}_{[\cdot,I]}  )^{-1}$ \cite{KailathEtAl2000}.

\subsection{Robust support-set based strategy}

The strategy described above performs a successive cancellation of
reconstructed signal coefficients. The performance is therefore
crucially dependent on detecting an active coefficient individually at
each stage, which is a `hard' decision. But the hypothesis of one
remaining active coefficient at each stage induces a risk of
irreversible detection errors. This increases with more erratically
evolving sparsity patterns, since the process is harder to predict.

The signal to prediction error ratio $\rho_i$, however, provides a statistic that can be viewed as `soft
information'. In order to increase robustness to detection errors we
propose to use $\rho_i$ for selecting the $\ell$ most likely remaining
coefficients. Let us denote the set of $\ell$ most likely indices by
$L$. It is joined with the existing set $I$ to form a hypothesized
support set $\tilde{I} = I\cup L$. This set is used to reconstruct
$\check{\mathbf{x}}_{[\tilde{I}]}$ and the coefficient $\check{x}_i$, $i
\in L$ with maximum magnitude is added to the support set $I$ at each
stage. Here $\ell = \max(0, K_{\text{max}} - |I|)$, which prevents
overfitting beyond the prior knowledge of the sparsity level.

Algorithm~\ref{alg:alternative} describes this
alternative detection strategy, based on a hypothesized support
set. The concerned scheme is called `robust predictive OMP' (rPrOMP).
\begin{algorithm}
\caption{: Robust predictive OMP (rPrOMP)}
\label{alg:alternative}
\begin{algorithmic}[1]
 \State Given: $\mathbf{y}, \mathbf{H}, \mathbf{R}^{-1},
 \hat{\mathbf{x}}^-$ and $\mathbf{P}^-$
   \State Set $k=0$, $\mathbf{r}_0 = \mathbf{y}_t$, $I = \varnothing$
    and $\mathbf{D}^{-1} = \mathbf{R}^{-1}$
    \Repeat
    \State $k := k+1$
    \State Compute $\rho_i$ using \eqref{eq:SPER} and \eqref{eq:SPERapprox}
    \State $\ell = \max(0, K_{\text{max}} - |I|)$
    \State $L = \{$indices of $\ell$ largest $ \rho_i \in I^c \}$
    \State $\tilde{I} = I \cup L$
   \State $\check{\mathbf{x}}_{[\tilde{I}]}
   =$\texttt{mmse-rec}$(\mathbf{y},\mathbf{H},\mathbf{R}^{-1},\hat{\mathbf{x}}^-,\mathbf{P}^-,\tilde{I})$
   \State $i_k = \argmax_{i \in L} | \check{x}_{i} |$
    \State $I := I \cup i_k$
    \State $\hat{\mathbf{x}}_{[I]}
    =$\texttt{mmse-rec}$(\mathbf{y},\mathbf{H},\mathbf{R}^{-1},\hat{\mathbf{x}}^-,\mathbf{P}^-,I)$;
    $\hat{\mathbf{x}}_{[I^c]} = \mathbf{0}$
    \State $\mathbf{r}_{k} = \mathbf{y} - \mathbf{H}_{[\cdot,I]}
    \hat{\mathbf{x}}_{[I]}$
   \State $\mathbf{D}^{-1} =$\texttt{update-cov}($\mathbf{D}^{-1},
   \mathbf{H},\mathbf{R}^{-1}, \mathbf{P}^-,I, i_k$)
    \Until{$(\| \mathbf{r}_{k} \|_2 \geq \| \mathbf{r}_{k-1} \|_2)$ or $(k
      > K_{\text{max}})$}
    \State Output: $\hat{\mathbf{x}}$ and $I$
\end{algorithmic}
\end{algorithm}

\subsection{Prediction of dynamic sparse signals}

Suppose a snapshot $\mathbf{y}_t$ has been observed and a prediction
$\hat{\mathbf{x}}^-_t$ is given along with $\mathbf{P}^-_t$. Let
$\hat{\mathbf{x}}_t$ denote the estimated sparse state vector after
the application of a predictive greedy pursuit algorithm (either PrOMP or rPrOMP), and $I$ its
support set. Then the updated error covariance matrix
$\mathbf{P}_t$ is computed block-wise corresponding to the set $I$ and its complement $I^c$. First,
$\mathbf{P}_{[I,I],t} = \left( \mathbf{S}^{-1}_t + \mathbf{H}^*_{[\cdot,I]} \mathbf{R}^{-1}_t \mathbf{H}_{[\cdot,I]} \right)^{-1}$ is the posterior error covariance, where $\mathbf{S}_t = \mathbf{P}^-_{[I,I],t}$ \cite{KailathEtAl2000}.
Second, the uncertainty of the inactive coefficients is preserved by
$\mathbf{P}_{[I^c,I^c],t} = \mathbf{P}^-_{[I^c,I^c],t}$. Finally, in line with
the MMSE reconstruction \eqref{eq:MMSE}, the
cross-correlations are set as $\mathbf{P}_{[I,I^c],t} =
\mathbf{0}$ and $\mathbf{P}_{[I^c,I],t} = \mathbf{0}$.

We propose predicting $\mathbf{x}_{t+1}$ from $\hat{\mathbf{x}}_t$ as a superposition of all possible transitions,
\begin{equation}
\hat{x}^-_{i,t+1} = \sum^{N}_{j=1} \lambda_{ji} \alpha_{ij} \hat{x}_{j,t},
\end{equation}
or written compactly, $\hat{\mathbf{x}}^-_{t+1} = \mathbf{F}
\hat{\mathbf{x}}_{t}$, where $f_{ij} = \lambda_{ji} \alpha_{ij}$. The
prediction error covariance matrix is then approximated by the equation,
$\mathbf{P}^-_{t+1} = \mathbf{F} \mathbf{P}_t \mathbf{F}^* +
\mathbf{Q}$.

Putting these blocks together we develop a Kalman-filter based algorithm
for recovery of sparse processes in
Algorithm~\ref{alg:SequentialPrOMP}, which we call dynamic
iterative pursuit (DIP). In DIP we use predictive OMP (PrOMP). If
robust predictive OMP (rPrOMP) is used instead, the algorithm can be referred to as `rDIP'.
\begin{algorithm}
\caption{: Dynamic Iterative Pursuit (DIP)}
\label{alg:SequentialPrOMP}
\begin{algorithmic}[1]
    \State Initialization $\hat{\mathbf{x}}^-_{0}$ and $\mathbf{P}^-_{0}$ \For{$t= 0, \dots$}
    \State \emph{\%Measurement update}
    \State $[  \mathbf{x}_t, I ]=$\texttt{PrOMP}$(\mathbf{y}_t, \mathbf{H}, \mathbf{R}^{-1}_t, \hat{\mathbf{x}}^-_t,\mathbf{P}^-_t)$
    \State $\mathbf{S}_t = \mathbf{P}^-_{[I,I],t}$
    \State $\mathbf{P}_{[I,I],t} = \left( \mathbf{S}^{-1}_t + \mathbf{H}^*_{[\cdot,I]} \mathbf{R}^{-1}_t \mathbf{H}_{[\cdot,I]} \right)^{-1}$
    \State  $\mathbf{P}_{[I^c,I^c],t} = \mathbf{P}^-_{[I^c,I^c],t}$; $\mathbf{P}_{[I,I^c],t} = \mathbf{0}$; $\mathbf{P}_{[I^c,I],t} = \mathbf{0}$
    \State \emph{\%Prediction}
    \State $\hat{\mathbf{x}}^-_{t+1} = \mathbf{F} \hat{\mathbf{x}}_{t}$
    \State $\mathbf{P}^-_{t+1} = \mathbf{F}\mathbf{P}_{t}
    \mathbf{F}^*  + \mathbf{Q}$
\EndFor
\end{algorithmic}
\end{algorithm}

\section{Experiments and Results}

In this section we evaluate DIP with respect to static OMP, SP and
convex relaxation based basis pursuit denoising (BPDN)
\cite{Candes&Wakin2008} algorithms. We also show the performance
of a `genie-aided' Kalman filter (KF) which provides a bound for
MMSE-based reconstruction of linear processes. The genie-aided
approach is given the sparsity pattern a priori, but does not know
the active signal coefficients. Finally, the robustness
properties of rDIP are compared with DIP for erratically
evolving sparsity patterns. The results are shown using Monte
Carlo simulations, averaged over 100 runs.

\subsection{Signal generation and performance measure}

Using a typical setup we consider a sparse process with the parameters $N =200$, $K=10$ and
number of snapshots $T=100$, with oscillating coefficients
according to an AR-model as in \eqref{eq:linearjump} with
$\alpha_{ij} = \alpha \equiv -0.8$, and $\mathbf{Q} = \sigma^2_w
\mathbf{I}_N$. The sparsity pattern transitions, $I_{x,t}
\rightarrow I_{x,t+1}$, are determined by transition probabilities
$\lambda_{ji}$ which are set in the experiments.

The transition of each active state $j \in I_{x,t}$ is generated
by a first-order Markov chain with $\lambda_{ji}$. If two states
in $I_{x,t}$ happen to transition into one, a new state is randomly
assigned to $I_{x,t+1}$, to ensure that the sparsity level is constant
in the experiment.

The entries of the sensing matrix $\mathbf{H}$ are set by
random drawing from a Gaussian distribution $\mathcal{N}(0,1)$
followed by unit-norm column scaling. The measurement noise
covariance matrix has form $\mathbf{R}_t = \sigma^2_n
\mathbf{I}_M$. Process and measurement noise are generated as
$\mathbf{w}_t \sim \mathcal{N}(\mathbf{0},\mathbf{Q})$ and
$\mathbf{n}_t \sim \mathcal{N}(\mathbf{0},\mathbf{R}_t)$,
respectively.

In the experiments, two signal parameters are varied; (a) the
signal-to-measurement noise ratio,
\begin{equation}
\text{SMNR} \triangleq \frac{\E \left[ \sum_t \| \mathbf{x}_t \|^2_2 \right]
}{\E \left[ \sum_t \| \mathbf{n}_t \|^2_2 \right] },
\end{equation}
while fixing $\E[ \| \mathbf{x}_t \|^2_2 ] \equiv 1$ so that
$\sigma^2_n = \frac{1}{M \times \text{SMNR}}$, and (b) the
fraction of measurements $\kappa \triangleq M/N$.

For a performance measure we use the signal-to-reconstruction
error ratio, defined as
\begin{equation}
\text{SRER} \triangleq \frac{\E \left[ \sum_t \| \mathbf{x}_t
\|^2_2
  \right] }{\E \left[ \sum_t  \| \mathbf{x}_t - \hat{\mathbf{x}}_t
    \|^2_2 \right] },
\end{equation}
which is the inverse of the normalized MSE. Note that SRER = 0~dB,
i.e. no reconstruction gain, is equivalent to using
$\hat{\mathbf{x}}_t = \mathbf{0}$.

\subsection{Algorithm initialization}

For the predictive algorithms---DIP, rDIP and
genie-aided KF---we use the mean and variance of an autoregressive
process as initial values, $\hat{\mathbf{x}}^-_{0} = \mathbf{0}$
and $\mathbf{P}^-_{0} = \sigma^2_x \mathbf{I}_N$ where $\sigma^2_x
= \frac{\sigma^2_w}{1-\alpha^2}$. In these algorithms we set
$K_{\text{max}} = K$ for consistent comparisons, although strict
equality is not a requirement.

Here we mention that BPDN \cite{Candes&Wakin2008} solves
\begin{equation*}
\hat{\mathbf{x}}_t = \argmin_{\mathbf{x}_t \in \mathbb{R}^N} \|
\mathbf{x}_t \|_1 \text{ subject to }  \| \mathbf{y}_t -
\mathbf{H}\mathbf{x}_t \|_2 \leq \varepsilon,
\end{equation*}
where the slack parameter $\varepsilon$ is determined by the
measurement noise power, as
\begin{equation*}
\varepsilon = \sqrt{\sigma^2_n (M+(2 \sqrt{2M}))},
\end{equation*}
following \cite{CandesEtAl2006,ChatterjeeEtAl2011a}. Note that
BPDN does not provide a $K$-element solution. It is also unable to
use prediction. The code for BPDN is taken from the $l_1$-magic
toolbox.

\subsection{Results}

For all experiments we ran 100 Monte Carlo simulations, where a new realization of $\{\mathbf{x}_t,\mathbf{y}_t \}^T_{t=1}$ and $\mathbf{H}$ was generated for each run.

\begin{figure}
  \begin{center}
    \includegraphics[width=1.00\columnwidth]{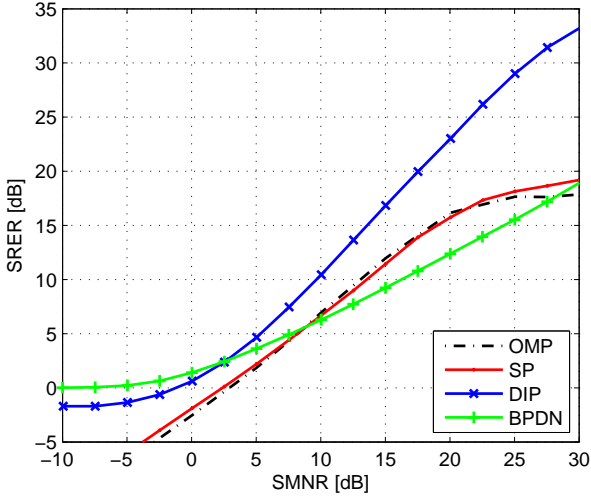}
  \end{center}
  \vspace*{-0.4cm}
  \caption{Comparison of different methods where we show SRER versus SMNR at $\kappa = 0.25$ and set transition probabilities
   according to \eqref{eq:lambda_model_1}.}
  \label{fig:SMNR}
\end{figure}

\begin{figure}
  \begin{center}
    \includegraphics[width=1.00\columnwidth]{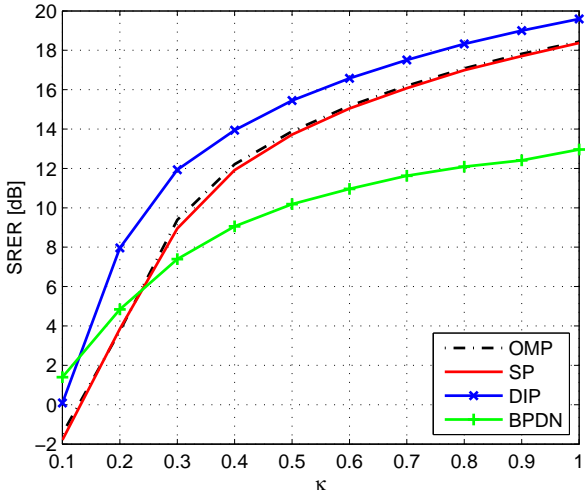}
  \end{center}
  \vspace*{-0.4cm}
  \caption{SRER versus $\kappa = M/N$ at SMNR = 10~dB. Transition
    probabilities according to \eqref{eq:lambda_model_1}.}
  \label{fig:kappa_omp}
\end{figure}

\begin{figure}
  \begin{center}
    \includegraphics[width=1.00\columnwidth]{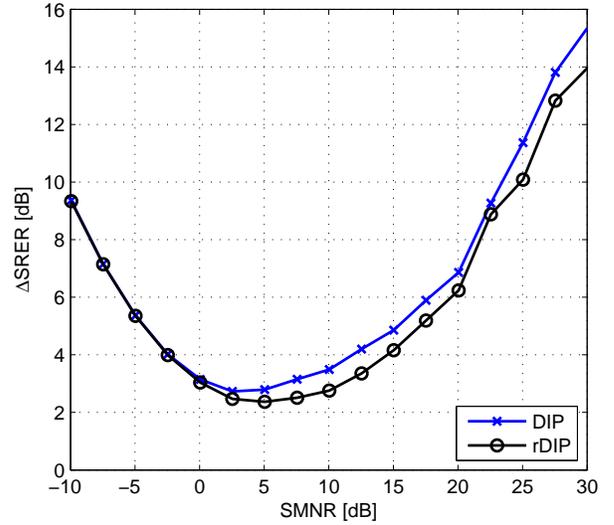}
  \end{center}
  \vspace*{-0.4cm}
  \caption{SRER improvement over OMP versus SMNR at $\kappa = 0.25$. Transition
    probabilities according to \eqref{eq:lambda_model_1}.}
  \label{fig:SMNR_gain}
\end{figure}

In the first experiment we consider a slowly varying sparsity pattern $I_{x,t}$
following the transition probabilities $\lambda_{ji}$ in \eqref{eq:lambda_model_1}. Figure \ref{fig:SMNR} shows how the algorithms perform with varying measurement noise
power at a fixed fraction of measurements $\kappa=0.25$. DIP overtakes
static BPDN at lower SMNR levels, while exhibiting a similar graceful
degradation. The static OMP and SP do not take into account the
measurement noise and hence continue to degrade. For instance, DIP
reaches the cut-off point of 0~dB reconstruction gain at an SMNR level
that is approximately 5~dB lower than the static
OMP. Figure~\ref{fig:kappa_omp} shows how the improvements persist for
varying $\kappa$ at a fixed SMNR = 10~dB.
In this scenario rDIP exhibits similar performance as DIP. Taking
static OMP as the baseline algorithm, the improvement of predictive
iterative pursuit is illustrated in Figure~\ref{fig:SMNR_gain}. The minimum SRER advantage is about 2~dB, and increases substantially with rising SMNR.

\begin{figure}
  \begin{center}
    \includegraphics[width=1.00\columnwidth]{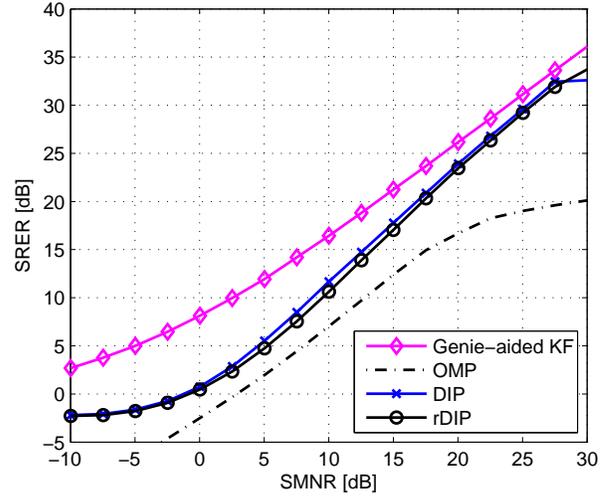}
  \end{center}
  \vspace*{-0.4cm}
  \caption{SRER versus SMNR at $\kappa = 0.25$. Degenerate transition
    probabilities (fixed sparsity pattern).}
  \label{fig:SMNR_kf}
\end{figure}

Next, we consider an unknown but static sparsity pattern $I_{x,t} \equiv I_x$, generated by degenerate transition probabilities $\lambda_{ji} = \delta(i-j)$, and compare DIP with a `genie-aided' KF. The latter filters the coefficients of a known support set $I_{x,t}$, and therefore provides an upper bound on the performance of sequential estimation. The bound is not necessarily tight since only $|I_{x,t}| \leq K$ is given in the problem. As
SMNR increases, DIP rapidly approaches the bound while OMP saturates
for $\kappa = 0.25$, illustrated in Figure~\ref{fig:SMNR_kf}. At SMNR
= 20~dB, OMP and DIP are about 10 and 2~dB from the upper limit,
respectively. Again, rDIP performs similarly to DIP.

\begin{figure}
  \begin{center}
    \includegraphics[width=1.00\columnwidth]{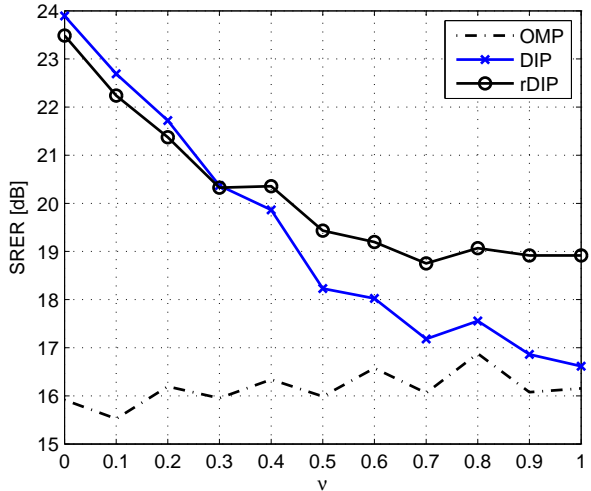}
  \end{center}
  \vspace*{-0.4cm}
  \caption{SRER versus mixture factor $\nu$ at SMNR = 20~dB
    and $\kappa = 0.25$. Transition probabilities according to
    \eqref{eq:lambda_model_2}.}
  \label{fig:nu_OMP}
\end{figure}

Finally, we consider an erratically evolving sparsity pattern
$I_{x,t}$, with transition probabilities $\lambda_{ji}$ set according to
\eqref{eq:lambda_model_2}, in order to compare the robustness of rDIP
with DIP. Figure~\ref{fig:nu_OMP} shows how performance is affected as
the mixture factor $\nu$ increases. DIP converges to OMP from above;
rDIP provides near equivalent performance to DIP at first but shows a
more graceful degradation. At the extreme, when all transitions are
equiprobable, rDIP is still capable of yielding above +2.5~dB
gain over OMP. This validates the robustness considerations behind
its design.

\emph{Reproducible results}: \textsc{Matlab} code for the
algorithms can be downloaded from
http://sites.google.com/site/saikatchatt/softwares. The codes
produce the results in Figures \ref{fig:SMNR} and \ref{fig:nu_OMP}.

\section{Conclusions}
We have developed a new iterative pursuit algorithm that uses sequential
predictions for dynamic compressive sensing, which we call dynamic
iterative pursuit. It incorporates prior statistical information
using linear MMSE reconstruction and the signal to prediction error as
a statistic. The algorithm was experimentally tested on a sparse signal with
oscillating coefficients and evolving sparsity pattern. The
results show that the algorithms exhibit graceful degradation at
low SMNR regions while capable of yielding substantial performance
gains as the SMNR level increases.

\section{Acknowledgement}
The authors would like to thank E. Candes and J. Romberg making
the $l_1$-magic toolbox available online.

\bibliography{refs_compress}
\bibliographystyle{ieeetr}

\end{document}